\documentclass[11pt,reqno]{amsart}

\newcommand{\mysection}[1]{\section{#1}
      \setcounter{equation}{0}}

\newtheorem{theorem}{Theorem}[section]
\newtheorem{lemma}[theorem]{Lemma}

\newtheorem{corollary}[theorem]{Corollary}

\theoremstyle{definition}

\newtheorem{definition}{Definition}[section]
\theoremstyle{remark}
\newtheorem{remark}{Remark}[section]

\newcommand\cbrk{\text{$]$\kern-.15em$]$}} 
\newcommand\opar{\text{\raise.2ex\hbox{${\scriptstyle | }$}\kern-.34em$($} }
\newcommand{\tr}{\text{\rm tr}\,}

 \makeatletter
 \def\dashint{%
 \operatorname%
 {\,\,\text{\bf--}\kern-.98em\DOTSI\intop\ilimits@\!\!}}
 \makeatother

\newcommand\bR{\mathbb{R}}

\newcommand\bX{\mathbb{X}}

\newcommand\cB{\mathcal{B}}
\newcommand\cF{\mathcal{F}}

\newcommand\cP{\mathcal{P}}

\newcommand\cD{\mathcal{D}}

\newcommand\frD{\mathfrak{D}}

\begin{document}

\title[It\^o-Wentzell formula]{
On  the It\^o-Wentzell formula for
distribution-valued processes and related topics}

\author[N.  Krylov]{N.V. Krylov}%
\thanks{The work   was partially supported
  by NSF grant DMS-0653121}
\address{127 Vincent Hall, University of Minnesota,
Minneapolis,
       MN, 55455, USA}
\email{krylov@math.umn.edu}


\begin{abstract}
We prove the It\^o-Wentzell formula for
 processes with values in the space of generalized functions
 by using
the stochastic Fubini theorem and
the It\^o-Wentzell formula for
real-valued processes, appropriate versions of which are also proved.
\end{abstract}

\maketitle

\mysection{Introduction and main result}

Let $(\Omega,\cF,P)$ be a complete probability space
with an increasing filtration $\{\cF_{t},t\geq0\}$
of complete with respect to $(\cF,P)$ $\sigma$-fields
$\cF_{t}\subset\cF$. Denote by $\cP$ the predictable
$\sigma$-field in $\Omega\times(0,\infty)$
associated with $\{\cF_{t}\}$ and let $\tau$ be
a stopping time with respect to $\{\cF_{t},t\geq0\}$. Let
 $w^{k}_{t}$, $k=1,2,...$, be independent one-dimensional
Wiener processes with respect to $\{\cF_{t}\}$.  
Let $\cD$ be the space 
of
generalized functions on the Euclidean 
$d$-dimensional space $\bR^{d}$
of points $x=(x^{1},...,x^{d})$. 

The following are just
versions of Definitions 4.4 and 4.6 of \cite{Kr99}.
Set $\bR_{+}=[0,\infty)$. Recall that for any $v\in\cD$
and $\phi\in C^{\infty}_{0}=C^{\infty}_{0}(\bR^{d})$ the function $(v,\phi(\cdot-x))$
is infinitely differentiable with respect to $x$, so that 
the sup in \eqref{11.16.2} below is predictable.

\begin{definition} 
 Denote by $\frD$                   \label{def 10.25.1}
 the set of all $\cD$-valued 
functions $u$ (written 
as $u_{t}(x)$ in a common abuse of notation)
on $\Omega\times\bR_{+}$ such that, for any $\phi\in 
C_{0}^{\infty}$,
 the restriction of the function $(u_{t},\phi)$ 
on $\Omega\times(0,\infty)$ is $\cP$-measurable 
and $(u_{0},\phi)$ is $\cF_{0}$-measurable.
 For $p=1,2$
denote by $\mathfrak{D}^{p}$ the subset of $\frD$
consisting of $u$ such that  
  for any  
  $\phi\in C_{0}^{\infty}$  and
 $T ,R \in\bR_{+}$, we have
\begin{equation}
                                            \label{11.16.2}
\int_{0}^{T}\sup_{ |x|\leq R}|(u_{t} ,
\phi(\cdot-x))|^{ p}\,dt<\infty
\quad\hbox{ (a.s.)}.
\end{equation}
In the same way, considering $\ell_{2}$-valued
distributions $g$  on $C_{0}^{\infty}$, that is
linear $\ell_{2}$-valued functionals
such that $(g,\phi)$ is continuous as an $\ell_{2}$-valued
function with respect to the standard convergence of
test functions, we define 
$\frD(\ell_{2})$
 and  $\mathfrak{D}^{ 2}
(\ell_{2})$ 
replacing
$|\cdot|$ in (\ref{11.16.2}) 
with $p=2$  by $|\cdot|_{\ell_{2}}$.

\end{definition}

Observe that if $g\in\mathfrak{D}^{ 2}(l_{2})$,
 then
for any $\phi\in C_{0}^{\infty}$,   and $T\in\bR_{+}$
$$
\sum_{k=1}^{\infty}\int_{0}^{T}(g^{k}_{t},\phi)^{2}\,dt
=\int_{0}^{T}|(g _{t},\phi)|_{\ell_{2}}^{2}\,dt<\infty
\quad\hbox{ (a.s.)},
$$
which, by well known theorems about convergence of
series of martingales, implies that the series in \eqref{12.23.40} below
converges uniformly on $[0,T]$ in probability for any $T\in\bR_{+}$.

\begin{definition} 
                                           \label{def 10.25.3}
Let $f, u\in\mathfrak{D}$, 
$g\in\mathfrak{D} (l_{2})$.
 We say that the equality
\begin{equation}
                                           \label{11.16.3}
du_{t}(x)=f_{t}( x)\,dt+
g_{t}^{k}( x)\,dw^{k}_{t},\quad t\leq\tau,
\end{equation}
holds {\em in the sense of distributions\/} if  
$ fI_{\opar0,\tau\cbrk}\in\mathfrak{D}^ {1}$, 
$gI_{\opar0,\tau\cbrk}\in\mathfrak{D}^{ 2}(l_{2})$ and 
for
 any $\phi\in C_{0}^{\infty}$, 
 with probability one  we have for all $t\in\bR_{+}$
\begin{equation}
                                             \label{12.23.40}
(u_{t\wedge\tau} ,\phi)=(u_{0} ,\phi)+\int_{0}^{t}I_{
s\leq\tau}
(f_{s},\phi)\,ds+\sum_{k=1}^{\infty}
\int_{0}^{t}I_{
s\leq\tau}(g^{k}_{s},\phi)\,dw^{k}_{s}.
\end{equation}
\end{definition}

Let $x_{t}$ be an $\bR^{d}$-valued stochastic process given by
$$
x^{i}_{t}=\int_{0}^{t}b^{i}_{s}\,ds
+\sum_{k=1}^{\infty}\int_{0}^{t}\sigma^{ik}_{s}\,dw^{k}_{s},
$$
where $b_{t}=(b^{i}_{t}),\sigma^{k}_{t}
=(\sigma^{ik}_{t})$ are predictable $\bR^{d}$-valued
processes such that for all $\omega$ and $ 
s, T\in\bR_{+}$
 we have  $\tr a_{s}<\infty$ and 
\begin{equation}
                                             \label{3.30.7}
\int_{0}^{T}(|b_{t}|+\tr a_{t})
\,dt<\infty,
\end{equation}
where $a_{t}=(a^{ij}_{t})$ and $2a^{ij}_{t}=(\sigma^{i\cdot},
\sigma^{j\cdot})_{\ell_{2}}$,
so that
$$
2\tr a_{t}=\sum_{i=1}^{d}\sum_{k=1}^{\infty}|\sigma^{ik}_{t}|^{2}.
$$

Finally, before stating our main result
we remind  the reader that for a generalized function $v$,
and any $\phi\in C^{\infty}_{0} $ the function
$(v,\phi(\cdot-x))$ is infinitely differentiable 
and for any derivative operator $D $ of order $n$ with respect 
to $x$ we have 
\begin{equation}
                                                 \label{4.3.2}
D (v,\phi(\cdot-x)) =(-1)^{n}(v,(D \phi)(\cdot-x))=:(D v,\phi(\cdot-x))=:
((D v)(\cdot+x),\phi)
\end{equation}
implying that $D u \in\frD $ if $u \in\frD $.

Here is our main result, that is a version 
of Lemma 4.7 of \cite{Kr99}.
 In case  $b\equiv0$ a proof of this lemma  is provided 
 in \cite{Kr99} 
without giving any precise indication as to
which version of  the It\^o-Wentzell
formula is used. We will fill this gap here. Set
$$
D_{i}=\frac{\partial}{\partial x^{i}},\quad D_{ij}=D_{i}D_{j}.
$$
\begin{theorem}
                                   \label{theorem 11.16.5}
Let $ f,u\in\mathfrak{D}$, 
$g\in\mathfrak{D} (l_{2})$.
Introduce
$$
v_{t}(x)=u_{t}(x+x_{t})
$$
and assume that \eqref{11.16.3} holds
 (in the sense of distributions). Then
$$
dv_{t}(x)= [f_{t}( x+x_{t})+a^{ij}_{t}
D_{ij}v_{t}( x)+b^{i}_{t}D_{i}v_{t}(x)
+(D_{i}g_{t}( x+x_{t}) ,\sigma^{i\cdot}_{t} )_{\ell_{2}}]\,dt
$$
\begin{equation}
                                                  \label{4.4.5}
+
[g^{k}_{t}( x+x_{t})+D_{i}v_{t}( x)\sigma^{ik}_{t} ]\,dw_{t}^{k},
\quad t\leq\tau
\end{equation}
(in the sense of distributions).
\end{theorem}

The reader understands that the summation convention
over the repeated indices $i,j=1,...,d$ 
(and $k=1,2,...$) is enforced here and throughout the
article.
The fact that \eqref{4.4.5} makes sense and indeed holds is proved in
Section \ref{section 4.6.1}. Our proof is outlined in \cite{Kr99}
and is based on the stochastic Fubini theorem and the It\^o-Wentzell
formula for real-valued processes. We prove a version
of the stochastic Fubini theorem in Section \ref{section 4.8.1}.
The It\^o-Wentzell
formula for real-valued processes in the form we need is proved in
Section \ref{section 4.8.2}.

There is a quite extensive literature on the stochastic Fubini theorem
(see, for instance, \cite{Pr} and \cite{NV} and the references
therein). It is worth saying that with some effort by using estimates
like \eqref{3.28.1}
we could obtain our version of the theorem in a somewhat weaker form
from probably the first one given in \cite{KS} or from more sophisticated
versions in \cite{NV}. In this case we would work with stochastic
integrals depending on the parameter $x$ as in
\begin{equation}
                                                     \label{4.8.1}
\sum_{k=1}^{\infty}\int_{0}^{t}I_{
s\leq\tau}(g^{k}_{s},\phi(\cdot-x))\,dw^{k}_{s}
\end{equation}
and for each $t$ choose a jointly measurable
 function of $(\omega,x)$, which is equal to \eqref{4.8.1} (a.s.)
for almost any $x$. However, there is a much better modification
working for all $x$ and $t$, which in the case of one-dimensional
semimartingales  is described in the corollary
of Theorem IV.63 of \cite{Pr} and obtained by using a  method
introduced by Dol\'eans-Dade. This modification allows also
to investigate   the continuous dependence on $t$ of the integral of 
\eqref{4.8.1} with respect to $x$, which in the case of
one driving semimartingale is proved in Theorems IV.64 and IV.65 of \cite{Pr}.
Our basic tools are Theorem IV.63, its above mentioned corollary,
and Theorem IV.64 of \cite{Pr} and are much more elementary than
rather involved arguments in \cite{NV}, where the authors
treat a very general situation, which is not within the scope 
of the present article, by
  using $\gamma$-radonifying operators and
the fact that $L_1$-spaces possess the UMD$^-$ property.

We prove and use  the stochastic Fubini theorem only for functions
given on $\bR^{d}$ with Lebesgue measure. Its generalization
for arbitrary $\sigma$-finite measure spaces is straightforward,
and can be used, as in \cite{KS}, to transform 
conditional expectations of stochastic integrals. This comment
is appropriate, because, actually, for the purpose
of proving Theorem \ref{theorem 11.16.5} one does not need
our stochastic Fubini theorem since $(u,\phi(\cdot-x))$
is an infinitely differentiable function of $x$ and one could just
approximate the integrals with respect to $x$ by
  Riemann sums and then pass to the limit.
This would prove the integral form of \eqref{4.4.5}
as in \eqref{12.23.40} for each fixed $t$ (a.s.)
and then an  additional effort based on our Corollary
\ref{corollary 4.3.1} is still needed to show that
the integral form holds (a.s.) for all $t$ at once.

Passing to the discussion of the It\^o-Wentzell
formula for real-valued processes  notice that 
our Theorem \ref{theorem 3.30.1} is somewhat close
to Theorem 3.3.1 of \cite{Ku}, 
which requires  two  derivatives of $F_{t}(x)$  in
$x$ to be continuous in $(t,x)$. Even if $F_{t}(x)$
is nonrandom, when the It\^o-Wentzell formula
becomes just It\^o's formula, our result is more general
than standard versions of It\^o's formula. For instance,
at those instances of time when $a_{t}=0$ we do not need
the second derivative of $F_{t}(x)$ to exist. 

Finally, it is worth pointing out that our results
are also true when there is only finitely many Wiener processes.
Considering infinitely many of them becomes
indispensable in the applications of the theory
of SPDEs to super-diffusions (see, for instance,
\cite{Kr97}).

\mysection{A version of the stochastic Fubini theorem}
                                                 \label{section 4.8.1}

If $E$ is a Borel subset   of a Euclidean space, by $\cB(E)$ we denote
the $\sigma$-field of Borel subsets of $E$. 
  Let $\Gamma$ be a Borel subset
of $\bR^{d}$ with nonzero finite Lebesgue measure.
\begin{definition}
                                              \label{definition 4.6.1}
Let
$B_{t}(x)$ be a real-valued function on
$\Omega\times\bR_{+}\times\Gamma $. We say that
it is a regular   field on $\Gamma$ if:

(a) It is
measurable with respect to
$\cF \otimes\cB(\bR_{+})\otimes\cB(\Gamma )$;

(b) For each $x\in\Gamma$, there is an event $\Omega_{x}$
such that  $P(\Omega_{x})=1$ and for any $\omega\in\Omega_{x}$, 
the function $B_{t}(\omega,x)$ is a continuous function of $t$ on
$\bR_{+}$; 

(c)  It is $\cF_{t}$-measurable for each $x\in\Gamma $
and $ t\in\bR_{+}$

We call it a regular martingale field on $\Gamma$ if
in addition

(d) For 
each $x\in\Gamma $
 the process $B_{t}(x)$ is a  local 
$\cF_{t}$-martingale on $\bR_{+}$ starting at zero.

\end{definition}

\begin{lemma}
                                           \label{lemma 4.12.1} 
If $B_{t}(x)$ is a regular field on $\Gamma$, then
there exists a regular field $A_{t}(x)$ on $\Gamma$ such that,
for each $x$, with probability one $A_{t}(x)=B_{t}(x)$
for all $t$ and

(b$'$) For each $\omega\in\Omega$ and $x\in\Gamma$ the function
$A_{t}(x)$ is continuous on $\bR_{+}$.

\end{lemma}

Proof.
By considering the processes $B^{n}_{t}(x)$, $n=1,2,...$, which are defined
as $(k+1-nt)B_{k/n}(x)+(nt-k)B_{(k+1)/n} (x)$
 for $k\leq nt\leq k+1$,
$k=0,1,...$,   noticing that, for each $x\in\Gamma$,
 $B^{n}_{t}(x)\to B_{t}(x)$ uniformly 
on each finite time interval in probability, and using
Theorem IV.62 of \cite{Pr} one easily obtains
a function $A_{t}(x)$ possessing the properties (a),   (b$'$),
and such that for each $x$, with probability one $A_{t}(x)=B_{t}(x)$
for all $t$. The latter and the completeness
 of $\cF_{t}$ implies that $A_{t}(x)$ also possesses property
(c). The lemma is proved.
\begin{definition}
                                                  \label{definition 4.12.1}
If a regular field on $\Gamma$ possesses property (b$'$)
of Lemma \ref{lemma 4.12.1}, then we call it strongly
regular.

\end{definition}
\begin{lemma}
                                           \label{lemma 3.28.1} 

Let  $p\in(0,\infty)$ and let
$m_{t}(x)$ be a regular martingale field on $\Gamma$.
Then
there exists a nonnegative strongly regular field $A_{t}(x)$ on
$\Gamma $ such that, for
each $x\in\Gamma $,
  with probability one $A_{t}(x)=\langle m(x)\rangle
_{t}$ for all $t\in\bR_{+}$. 

Moreover, if
$A_{t}(x)$ is a function with the above described
properties and such that

(i)   It is $\cF_{t}\otimes\cB(\Gamma)$-measurable for 
each   $t\in \bR_{+}$;

(ii) Almost surely

\begin{equation}
                                             \label{3.28.3}
\int_{\Gamma }\sup_{ t \in\bR_{+} }
A_{t}^{p/2}(x)\,dx<\infty ,
\end{equation}
then for any countable set $\rho\subset\bR_{+}$ with probability one
\begin{equation}
                                             \label{3.28.1}
\int_{\Gamma }\sup_{t\in\rho }|m_{t}(x)|^{p}\,dx<\infty
\end{equation}
and for any $\varepsilon,\delta>0$ we have
\begin{equation}
                                             \label{3.28.2}
P\big(\int_{\Gamma }\sup_{t\in\rho}|m_{t}(x)|^{p}\,dx
\geq\delta\big)\leq P(C_{\infty}\geq\varepsilon)+
 \frac{N}{\delta}
E(\varepsilon \wedge C_{\infty})
\end{equation}
where the constant $N$ depends only on $p$ and
$$
C_{t}:=\int_{\Gamma }\sup_{ s\leq t}A^{p/2}_{s}(x)\,dx.
$$
\end{lemma}

Proof. To prove the first assertion notice that by 
the corollary of Theorem IV.63  of \cite{Pr}
there exists a strongly regular martingale field
$B_{t}(x)$ on $\Gamma$ such that,
for each $x\in\Gamma $, with probability one
$$
B_{t}(x)=\int_{0}^{t}m_{t}(x)\,dm_{t}(x)
$$
for all $t$.  Actually, the corollary of Theorem IV.63 of \cite{Pr}
is stated somewhat differently, so that what we need
follows from its proof and the arguments leading to 
the corollary. Taking a strongly regular modification $n_{t}$ of 
$m^{2}_{t}(x)$, which exists by Lemma \ref{lemma 4.12.1} and
letting $A_{t}(x)=|n_{t}(x)-2B_{t}(x)|$
yields a function we are looking for.
 
To prove the second assertion, observe that 
the process $C_{t}$
is $\cF_{t}$-adapted and,
with probability one, is continuous in $t\in\bR_{+}$
owing to condition \eqref{3.28.3} and the dominated
convergence theorem.
Therefore
$$
\tau:= \inf\{t\geq0: C_{t}\geq\varepsilon\}
$$
is a stopping time. Now
$$
P\big(\int_{\Gamma }\sup_{t\in\rho}|m_{t}(x)|^{p}\,dx
\geq\delta\big)\leq P(C_{\infty}\geq \varepsilon) 
$$
$$
+P\big(\int_{\Gamma }\sup_{t\in\rho}|m_{t}(x)|^{p}\,dx
\geq\delta,\tau=\infty\big),
$$
where the last term is less than
$$
P\big(\int_{\Gamma }\sup_{t\in\rho}|m_{t\wedge\tau}(x)|^{p}\,dx
\geq\delta \big)\leq\frac{1}{\delta}
E\int_{\Gamma }\sup_{t\in\rho}|m_{t\wedge\tau}(x)|^{p}\,dx,
$$
which, in turn, by the Burkholder-Davis-Gundy
  inequalities is dominated by
$$
\frac{1}{\delta}
\int_{\Gamma }E\sup_{t }|m_{t\wedge\tau}(x)|^{p}\,dx
\leq\frac{N}{\delta}
\int_{\Gamma }EA^{p/2}_{\tau}(x)\,dx
\leq\frac{N}{\delta}
EC_{\tau}\leq\frac{N}{\delta}
E(\varepsilon \wedge C_{\infty}).
$$
This proves \eqref{3.28.2} which implies \eqref{3.28.1}
 if one first lets
$\delta\to\infty$ and then $\varepsilon\to\infty$.
The lemma is proved.

For any multi-index $\alpha=(\alpha_{1},...,\alpha_{d})$,
$\alpha_{1},...,\alpha_{d}\in\{0,1,...\}$, define
$$
D^{\alpha}=D_{1}^{\alpha_{1}}\cdot...\cdot D_{d}^{\alpha_{d}},
\quad |\alpha|=|\alpha_{1}|+...+|\alpha_{d}|.
$$
In the following corollary  
 $\Gamma$ is a ball, $p\in[1,\infty)$
 and  $n$ is an integer. We denote $\lambda=n-d/p$
and assume that either $p>1$ and $\lambda\in(0,1)$
or $p=1$ and $\lambda=1$, so that $n=d+1$.
\begin{corollary}
                                         \label{corollary 4.3.1}

(i) Let $m_{t}(x)$ be a regular martingale field on $\Gamma$ and  assume
that for each $\omega\in\Omega$ and $t\in\bR_{+}$
it is $n$ times continuously differentiable in~$x$.

(ii) Suppose that, for each multi-index $\alpha$ with $|\alpha|\leq
n$, $D^{\alpha}m_{t}(x)$ is also
a regular martingale field on $\Gamma$.

(iii)
Finally, assume that for each multi-index $\alpha$ 
with $|\alpha|\leq n$ (including $\alpha=0$) 
on $ \Gamma$
there
exists a nonnegative strongly regular field $A^{\alpha}_{t}(x)$
possessing the properties   (i)  and (ii) 
of Lemma \ref{lemma 3.28.1}
and such that, for each $x\in\Gamma$,
with probability one
 $A^{\alpha}_{t}(x)=\langle D^{\alpha}m(x)\rangle
_{t}$ for all $t\in\bR_{+}$. 

Then   there is a (finite) random
variable $\nu$ such that with probability one
for all $x,y\in\Gamma$ and $t\in\bR_{+}$ we have
\begin{equation}
                                              \label{4.3.1}
|m_{t}(x)-m_{t}(y)|\leq \nu|x-y|^{\lambda}.
\end{equation}
Furthermore, with probability one $m_{t}(x)$
is continuous with respect to $(t,x)$ on $\bR_{+}\times\Gamma$.
\end{corollary}

Proof. Take $\rho$ as the set of rational numbers
on $\bR_{+}$, and observe that, owing to \eqref{3.28.1},
there is an event $\Omega'$ of full probability and such
that for any $\omega\in\Omega'$   we have
$$
\sup_{t\in\rho}\sum_{|\alpha|\leq n}
\int_{\Gamma }|D^{\alpha}m_{t}(x)|^{p}\,dx=:\nu_{0}<\infty.
$$
By the Sobolev embedding
theorems (see, for instance, Theorem 5.4 of \cite{Ad}),
for each $\omega$ and $t$, for which
$$
\sum_{|\alpha|\leq n}
\int_{\Gamma }|D^{\alpha}m_{t}(x)|^{p}\,dx<\infty,
$$
there exists a continuous function $v(x)$ on $\Gamma$
such that $v(x)=m_{t}(x)$ for almost all $x\in\Gamma$
and
$$
|v(x)-v(y)|\leq N|x-y|^{\lambda}\sum_{|\alpha|\leq n}
\int_{\Gamma }|D^{\alpha}m_{t}(x)|^{p}\,dx\quad\forall x,y\in\Gamma,
$$
where $N$ depends only on $ d,p$, and $\Gamma $. Of course,
$v(x)=m_{t}(x)$ for  all $x\in\Gamma$, since $m_{t}(x)$
is assumed to be continuous in $x$. Therefore, \eqref{4.3.1}
holds with $\nu=N\nu_{0}$ for all rational $t$,
  $\omega\in\Omega'$, and $x,y\in\Gamma$.

Let $\bX$ be the set of points with rational coordinates in $\Gamma$
and for each $x\in\bX$ let $\Omega_{x}$ be the event of full probability
such that for each $\omega\in\Omega_{x}$ the function
$m_{t}(x)$ is continuous in $t$. Then for any 
$$
\omega\in\Omega''
:=\Omega'\bigcap_{x\in\bX}\Omega_{x}
$$
and $x,y\in\bX$
we have \eqref{4.3.1}  for all rational, and hence, for all $t$.
Since $m_{t}(x)$ is assumed to be continuous in $x$, in \eqref{4.3.1}
one can take arbitrary $x,y\in\Gamma$ and $t\in\bR_{+}$
 as long as $\omega\in\Omega''$. For those $\omega$ and any $x \in
\Gamma$ it holds that $m_{t}(x_{n})\to m_{t}(x )$ uniformly in $t$
if $x_{n}\to x $. By taking $x_{n}\in\bX$, so that $m_{t}(x_{n})$
are continuous in $t$, we conclude that $m_{t}(x )$ is continuous
in $t$ for any $\omega\in\Omega''$ and $x\in\Gamma$. Since
it is also uniformly continuous in $x$, it is jointly
continuous with respect to $(t,x)$ for $\omega\in\Omega''$.
It only remains to observe that
obviously $P(\Omega'')=1$ and this proves the corollary.
 
\begin{remark}
The above corollary is close in spirit to
Theorem 3.1.1 of \cite{Ku}. However, in the 
applications to the integrals like \eqref{4.8.1} we have in mind
(see, for instance,  Lemma  \ref{lemma 4.10.1})
it is much easier to use the corollary than
Theorem 3.1.1 of \cite{Ku}.
\end{remark}
\begin{corollary}
By taking $q\in(0,1)$, substituting $\delta^{1/q}$  
in \eqref{3.28.2} in place of $\delta$ and $\varepsilon$,
  and then integrating the result with respect to
$\delta$ over $(0,\infty)$, we obtain
$$
E\big(\int_{\Gamma }\sup_{t }|m_{t}(x)|^{p}\,dx\big)^{q}
\leq NE\big(\int_{\Gamma }\sup_{t }A^{p/2}_{t}(x) \,dx\big)^{q },
$$
where the constant $N$ depends only on $p$ and $q$.
\end{corollary}

Estimate \eqref{3.28.2}
allows us to improve in our particular case
Theorem 65 of \cite{Pr}, in which in condition \eqref{3.29.3} below
the power 1/2 is replaced with 1.

\begin{lemma}
                                           \label{lemma 3.29.1}
Let $f_{t}(x)$ be a real-valued function
on $\Omega\times(0,\infty)\times\Gamma $ which is
$\cP\otimes\cB(\Gamma )$-measurable and such that
$$
 \int_{0}^{\infty}f_{t}^{2}(x)\,dt<\infty 
$$
for each $x\in\Gamma $ and $\omega$. Then there exists  
a strongly regular martingale field $m_{t}(x)$  on
$ \Gamma $  such that for 
each $x\in\Gamma $ with probability one
$$
m_{t}(x)=\int_{0}^{t}f_{s}(x)\,dw_{s}
$$
for all $t$. Furthermore, if
\begin{equation}
                                             \label{3.29.3}
\int_{\Gamma }\big(\int_{0}^{\infty}
f^{2}_{t}(x)\,dt\big) ^{1/2}\,dx<\infty \quad\text{(a.s.)},
\end{equation}
then for any function $m_{t}(x)$ with the properties
described above
\begin{equation}
                                             \label{3.29.4}
\int_{0}^{\infty}\big(\int_{\Gamma }f_{s}(x)\,dx\big)^{2}
\,ds<\infty,\quad
\int_{\Gamma }\sup_{t}| m_{t}(x)| \,dx<\infty
 \quad\text{(a.s.)},
\end{equation}
the stochastic integral
\begin{equation}
                                             \label{3.29.5}
\int_{0}^{t}\big(\int_{\Gamma }f_{s}(x)\,dx\big)\,dw_{s}
\end{equation}
is well defined, and with probability one
\begin{equation}
                                             \label{3.29.2}
\int_{\Gamma } m_{t}(x) \,dx=
\int_{0}^{t}\big(\int_{\Gamma }f_{s}(x)\,dx\big)\,dw_{s}
\end{equation}
for all $t$.
\end{lemma}

Proof.
The existence of $m_{t}(x)$ with  the claimed properties
follows from the corollary of Theorem IV.63  of \cite{Pr}.
Furthermore, the function
$$
A_{t}(x)=\int_{0}^{t}f^{2}_{s}(x)\,ds
$$
is certainly a strongly regular field on $\Gamma$ 
such that, for any $x\in\Gamma$, with probability one
 $A_{t}(x)=\langle m(x)\rangle_{t}$
for all $t$. Furthermore, $A_{t}(x)$ possesses
property (i)   of Lemma 
\ref{lemma 3.28.1} and property (ii) with $p=1$ if condition
\eqref{3.29.3} is satisfied. Under this condition,
 which we assume in the rest of the proof,
the first inequality in \eqref{3.29.4}
follows from \eqref{3.29.3} by Minkowski's inequality
and implies that \eqref{3.29.5} is well defined indeed.
Also, \eqref{3.28.1} with $p=1$ yields
the second inequality in \eqref{3.29.4}.

Equality \eqref{3.29.2} follows from Theorem 64 of \cite{Pr}
if $f$ is bounded. 
In the general case for $n=1,2,...$ define
$\chi_{n}(s)=(-n)\vee s\wedge n$, $f^{n}_{t}=\chi_{n}
(f_{t})$, and let $m^{n}_{t}(x)$ be a strongly regular martingale field
on $\Gamma$
 such that for 
each $x\in\Gamma $ with probability one
$$
m^{n}_{t}(x)=\int_{0}^{t}f^{n}_{s}(x)\,dw_{s}
$$
for all $t$. By the above with probability one
\begin{equation}
                                             \label{3.29.7}
\int_{\Gamma } m^{n}_{t}(x) \,dx=
\int_{0}^{t}\big(\int_{\Gamma }f^{n}_{s}(x)\,dx\big)\,dw_{s}
\end{equation}
for all $t$. By \eqref{3.28.2} and the dominated convergence theorem
for any $\varepsilon,\delta>0$
$$
P\big(\int_{\Gamma } \sup_{t}|  m^{n}_{t}(x)-m_{t}(x)| \,dx
\geq \delta\big)
$$
$$
\leq P\big(\int_{\Gamma }\big(\int_{0}^{\infty}
| f^{n}_{t}(x)-f_{t}(x)|^{2}\,dt\big)^{1/2}\,dx\geq\varepsilon
\big)+N\varepsilon/\delta\to N\varepsilon/\delta
$$
as $n\to\infty$. Also by Minkowski's inequality
$$
\int_{0}^{\infty}\big(\int_{\Gamma}| f^{n}_{t}(x)-f_{t}(x)| \,dx
\big)^{2}\,dt\leq
\big(\int_{\Gamma}\big(\int_{0}^{\infty}
| f^{n}_{t}(x)-f_{t}(x)|^{2}\,dt\big)^{1/2}\,dx\big)^{2}
\to0
$$
as $n\to\infty$ for almost any $\omega$. Hence 
both sides
of \eqref{3.29.7} converge in probability to the corresponding
sides of \eqref{3.29.2} uniformly in $t$ implying that \eqref{3.29.2}
holds with probability one for all $t$ and the lemma is proved. 
\begin{remark}
                                                      \label{remark 4.6.1}

Below we are going to use    ``local" versions
of Lemmas \ref{lemma 3.28.1} and \ref{lemma 3.29.1}
when all processes will be considered on $[0,T]$ with a $T\in\bR_{+}$.
These versions are obtained by replacing $m_{t}(x)$ and $f_{t}(x)$
with $m_{t\wedge T}$ and $f_{t}I_{t<T}$ respectively.
 
\end{remark}

Here is a version of the stochastic Fubini theorem.

\begin{lemma}
                                          \label{lemma 3.30.2}
Let $T\in\bR_{+}$ and
let $ G_{t}(x)$ be real-valued  and $H_{t}(x)=(H^{k}_{t}(x),\\
k=1,2,...)$
be $\ell_{2}$-valued functions defined on $\Omega\times(0,T]
\times \Gamma$ and possessing the following properties:

(i) The functions $G_{t}(x)$ and $H_{t}(x)$ are $\cP_{T}\otimes
\cB(\Gamma)$-measurable, where $\cP_{T}$ is the restriction
of $\cP$ to $\Omega\times(0,T]$;

(ii)  There is an 
event $\Omega'$ of full probability
such that  for each   $\omega\in \Omega' $  and $x\in\Gamma$ we have
$$
\int_{0}^{T}(|G_{t}(x)|+|H_{t}(x)|^{2}_{\ell_{2}})\,dt<\infty;
$$

(iii) We have (a.s.)
$$
\int_{0}^{T}\int_{\Gamma} |G_{t}(x)|\,dxdt
+\int_{\Gamma}\big(\int_{0}^{T} |H_{t}(x)|^{2}_{\ell_{2}}
 \,dt\big)^{1/2}
\,dx <\infty.
$$

Under these assumptions we claim that

(a) There is a function $F_{t}(x)$ on $\Omega\times[0,T]\times\Gamma$,
which is $\cF \otimes\cB([0,T]) \otimes\cB(\Gamma)$-measurable,
continuous in $t$, and such that
for any $x\in\Gamma$ with probability one we have
\begin{equation}
                                              \label{4.1.1}
F_{t}(x)=\int_{0}^{t}G_{s}(x)\,ds+\sum_{k=1}^{\infty}\int_{0}^{t}
H^{k}_{s}(x)\,dw^{k}_{s}
\end{equation}
for all $t\in[0,T]$,
where the series converges uniformly on $[0,T]$ in probability;

(b) For any   $k=1,2,...$,
the stochastic integrals (no summstion in $k$)
$$
\int_{0}^{t}\int_{\Gamma}H^{k}_{s}(x)\,dxdw^{k}_{s}
$$
are well defined for $t\in[0,T]$;

(c) If we are given
a function $F_{t}(x)$ on $\Omega\times[0,T]\times\Gamma$
with somewhat weaker properties, namely,
such that

(iv) For each $t\in[0,T]$ the function $F_{t}(x)$ is
measurable in $(\omega,x)$ with respect to the completion 
$\overline{\cF\otimes\cB(\Gamma)}$ of
  $\cF \otimes\cB(\Gamma)$ with respect to the product
measure;

(v)   For each $t\in[0,T]$ and $x\in\Gamma$  
equation \eqref{4.1.1} holds almost surely,

\noindent then  for any countable subset  $\rho$ of $[0,T]$
\begin{equation}
                                                \label{3.30.5}
\int_{\Gamma}\sup_{ t\in\rho}|F_{t}(x)|\,dx<\infty
\quad\text{(a.s.)},
\end{equation}
  and 
for each $t\in[0,T]$  almost surely
\begin{equation}
                                                \label{3.30.6}
\int_{\Gamma} F_{t}(x)\,dx=\int_{0}^{t}\int_{\Gamma} G_{s}(x)\,dx
ds+\sum_{k=1}^{\infty}
\int_{0}^{t}\int_{\Gamma}H^{k}_{s}(x)\,dxdw^{k}_{s},
\end{equation}
  where the series converges uniformly
on $[0,T]$ in probability.  

(d) If for a function $F_{t}(x)$ as in (c),
 for almost all $(\omega,x)$, $F_{t}(x)$ is continuous in $t$
on $[0,T]$ (like the one from assertion (a)),
then with probability one \eqref{3.30.6} holds for all $t
\in[0,T]$.
\end{lemma}

Proof.   Obviously, replacing $G$ and $H$ with
$GI_{\Omega'}$ and $HI_{\Omega'}$, respectively,
 will not affect anything
and therefore we may assume that assumption (ii) holds with
$\Omega'=\Omega$.
The fact that for each $x$
 the series in \eqref{4.1.1}
converges uniformly on $[0,T]$
in probability due to condition (ii) is discussed after Definition
\ref{def 10.25.1}.
As there, the fact that, by Minkowski's
inequality
$$
\bigg(\sum_{k=1}^{\infty}\int_{0}^{T}\big(\int_{\Gamma}
H^{k}_{t}(x)\,dx\big)^{2}\,dt\bigg)^{1/2}
\leq \int_{\Gamma}\bigg(\sum_{k=1}^{\infty}\int_{0}^{T}|
H^{k}_{t}(x)|^{2}\,dt\bigg)^{1/2}\,dx,
$$
where the latter is finite (a.s.)  
due to (iii), implies that the series in \eqref{3.30.6}
converges uniformly
on $[0,T]$ in probability.

By Lemma \ref{lemma 3.29.1} (see also Remark \ref{remark 4.6.1})
for each $k$, on $\Omega\times
[0,T]\times\Gamma$ there exists an $\cF\otimes\cB([0,T])
\otimes\cB(\Gamma)$-measurable function $m^{k}_{t}(x)$,
which is
continuous in $t$ 
for each $x\in\Gamma $ and $\omega$ and
such that  for any $x\in\Gamma$ with probability one
$$
m^{k}_{t}(x)=\int_{0}^{t} H^{k}_{s}(x)\,dw^{k}_{s},
\quad
\int_{\Gamma}\sup_{s\leq t}|m^{k}_{s}(x)|\,dx<\infty
$$
for all $t\in[0,T]$. Furthermore,  
by Lemma \ref{lemma 3.29.1} we also have that with 
probability one
\begin{equation}
                                                  \label{4.1.2}
\int_{\Gamma}m^{k}_{t}(x)\,dx=\int_{0}^{t}
\int_{\Gamma}H^{k}_{s}(x)\,dx\,dw^{k}_{s}
\end{equation}
for all $t\in[0,T]$. Now introduce
$$
M^{k}_{t}(x)= \sum_{j=1}^{k}m^{j}_{t}(x).
$$
As we have pointed out in the beginning of the proof,
for each $x$,
 the processes
$M^{k}_{t}(x)$ converge uniformly on $[0,T]$
in probability as $k\to\infty$. By Theorem 62 of \cite{Pr}
there exists an $\cF\otimes\cB([0,T])
\otimes\cB(\Gamma)$-measurable function $m_{t}(x)$, which is
continuous in $t$ for all $\omega$ and $x$ and such that for any
$x\in\Gamma$ we have $M^{k}_{t}(x)\to m_{t}(x)$ uniformly on 
$[0,T]$ in probability as $k\to\infty$. 
Of course, for each $x\in\Gamma$ with probability one
$$
m_{t}(x)=\sum_{k=1}^{\infty}\int_{0}^{t}
H^{k}_{s}(x)\,dw^{k}_{s}
$$
for all $t\in[0,T]$ and this certainly
 proves   assertion (a). Assertion (b)
is proved above.

Next, condition (v) means that
for each $t\in[0,T]$ and $x\in\Gamma$ we have (a.s.)
\begin{equation}
                                                 \label{3.31.5}
F_{t}(x)=\int_{0}^{t}G_{s}(x)\,ds+m_{t}(x) .
\end{equation}

Furthermore,
for each $x$ the process $m_{t}(x)$ is a continuous
local martingale and
$$
\langle m(x)\rangle_{t}=\int_{0}^{t}|H_{s}(x)|_{\ell_{2}}^{2}\,ds
$$
for all $t\in[0,T]$ (a.s.). The right-hand side here 
can be  taken as $A_{t}(x)$ in Lemma \ref{lemma 3.28.1}
and this $A_{t}(x)$  
is strongly regular by the stipulation made in the
 beginning of the proof, 
satisfies condition (i) of that lemma and also satisfies
its condition (ii) with $p=1$
 due to condition (iii) of the present lemma.
By Lemma \ref{lemma 3.28.1} we have
\begin{equation}
                                                 \label{3.31.4}
\int_{\Gamma}\sup_{t\leq T}|m_{t}(x)|\,dx<\infty
\end{equation}
(a.s.). Furthermore, in light of \eqref{3.31.5}
for   each  $x$
\begin{equation}
                                                 \label{3.31.3}
\sup_{t\in\rho}|F_{t}(x)-\int_{0}^{t}G_{s}(x)\,ds-m_{t}(x)|=0 
\end{equation}
(a.s.).
Here the left-hand side is a $
\overline{\cF\otimes\cB(\Gamma)}$-measurable.
Therefore, for almost any $\omega$ equation \eqref{3.31.3} holds 
for almost all $x$. This, \eqref{3.31.4}, and condition (iii)
imply \eqref{3.30.5}.

Also, for the local martingale $m_{t}(x)-M^{k}_{t}(x)$
the process $\langle   m(x)-M^{k} (x)\rangle_{t}$
  can be taken to be
$$
A^{k}_{t}(x) =\sum_{j=k+1}^{\infty}
\int_{0}^{t}(H^{j}_{s}(x))^{2}\,ds,
$$
so that by Lemma \ref{lemma 3.28.1}
for any $\varepsilon,\delta>0$  
$$
P\big(\int_{\Gamma}\sup_{t 
\leq T}|m_{t}(x)-M^{k}_{t}(x)|\,dx
\geq\delta\big)\leq P\big(\int_{\Gamma}(A^{k}_{T}(x))^{1/2}\,dx
\geq\varepsilon\big)
 +\frac{N\varepsilon}{\delta}.
$$
After letting first $k\to\infty$ and then $\varepsilon
\to0$   we conclude from
assumption (iii) that
$$
\int_{\Gamma}\sup_{t\leq T }|m_{t}(x)-M^{k}_{t}(x)|\,dx\to0
$$
as $k\to\infty$ in probability. In particular,
for each $t\in[0,T]$
$$
\sum_{i=1}^{k}\int_{\Gamma} m^{k}_{t}(x) \,dx
=\int_{\Gamma} M^{k}_{t}(x) \,dx\to
\int_{\Gamma} m_{t}(x) \,dx
$$
as $k\to\infty$ in probability which is to say that (a.s.)
\begin{equation}
                                                 \label{3.31.7}
\sum_{i=1}^{\infty}\int_{\Gamma} m^{k}_{t}(x) \,dx
=\int_{\Gamma} m_{t}(x) \,dx.
\end{equation}
Now for $t\in[0,T]$ being fixed and for almost  any $\omega$
we have that equation
\eqref{3.31.5} holds for almost all $x$. We integrate
it over $\Gamma$ and using \eqref{3.31.7} and
\eqref{4.1.2} conclude that
\eqref{3.30.6} indeed holds (a.s.).
This finishes the proof of  assertion~(c).

Observe that the right-hand side of \eqref{3.30.6}
is continuous in $t$ on $[0,T]$ with probability one.
If for almost all $(\omega,x)$, $F_{t}(x)$ is continuous in $
[0,T]$,
then on the set of full probability for almost any 
$x\in\Gamma$ the function $F_{t}(x)$ is continuous on $[0,T]$
and owing to the dominated convergence theorem
and \eqref{3.30.5} the left-hand side of \eqref{3.30.6}
is also continuous in $t$ on $[0,T]$ with probability one.
This implies   assertion (d) of the lemma,
which is thus proved.

\mysection{A real-valued version of the It\^o-Wentzell formula}
                                                \label{section 4.8.2}

For $\gamma\in(0,\infty)$ set $B_{\gamma}=\{x\in\bR^{d}:|x|<\gamma\}$.
Also introduce
$$
L_{t}v=a^{ij}_{t}D_{ij} v+b^{i}_{t}D_{i}v,\quad
\Lambda^{k}_{t}v=\sigma^{ik}_{t}D_{i}v.  
$$

\begin{theorem}
                                          \label{theorem 3.30.1}
Let $T\in\bR_{+}$ and $\gamma\in(0,\infty)$. Set 
$\eta_{t}(x)=I_{B_{\gamma}}(x-x_{t})$.
Let  real-valued  $G_{t}(x)$  and
$\ell_2$-valued  $H_{t}(x)=(H^{k}_{t}(x),
k=1,2,...)$ be some 
functions on $\Omega\times(0,T]
\times\bR^{d}$
satisfying  assumption  (i)
of Lemma \ref{lemma 3.30.2} for any ball $\Gamma$
and let
real-valued $F_{t}(x)$ be a
function  on $\Omega\times[0,T]
\times\bR^{d}$   such that:

(i) For each $x$ the restriction of the function $F_{t}(x)$ 
to $\Omega\times(0,T]$ is 
$\cP_{T}$-measurable and $F_{0}(x)$
is $\cF_{0}$-measurable;

(ii) For   any $\omega\in\Omega$ and $t\in[0,T]$ 
 the function $F_{t}(x)$ is continuous in
$x$;

(iii) For almost any $(\omega,t)\in\Omega\times[0,T]$ and
$k=1,2,...$ 

\quad(a)
 the functions $G_{t}(x)$, $H^{k}_{t}(x)$, and
$|H_{t}(x)|_{\ell_{2}}$ 
are continuous in $x$,

\quad(b)   the generalized functions
$L_{t}F_{t}(x)$, $\Lambda^{k}_{t}F_{t}(x)$,
and $\Lambda^{k}_{t}H^{k}_{t}(x)$ are continuous functions of $x$
as well as the functions
$|\Lambda_{t}F_{t}(x)|_{\ell_{2}}$ and
$$
\sum_{r=1}^{\infty}|
\Lambda^{r}_{t}H^{ r}_{t}(x)|;
$$
 
(iv)  There is an event $\Omega'$
of full probability such that   for each $x\in\bR^{d}$
and $\omega\in \Omega' $
we have
$$
\int_{0}^{T}\big( \eta_{t}(x)|F_{t}(x)|(|b_{t}|+\tr a_{t})
+\eta_{t}(x) F_{t}^{2}(x) \tr a_{t}+|G_{t}(x)|+|H_{t}(x)|^{2}_{\ell_{2}}
\big)\,dt<\infty 
$$
 and, for each $x\in\bR^{d}$,
 with probability
one equation 
\begin{equation}
                                                           \label{4.9.1}
F_{t}(x)=
F_{0}(x)+
\int_{0}^{t}G_{s}(x)\,ds+\sum_{k=1}^{\infty}\int_{0}^{t}
H^{k}_{s}(x)\,dw^{k}_{s} 
\end{equation}
 holds for all $t\in[0,T]$;

(v) We have  
$$
\int_{\bR^{d}}\int_{0}^{T}  \eta_{t}(x)
 |F_{t}(x)|(|b_{t}|+\tr a_{t}) 
  \,dxdt
$$
\begin{equation}
                                                           \label{4.6.3}
+ \int_{\bR^{d}}\big(
\int_{0}^{T}\eta_{t}(x)
|F_{t}(x)|^{2}\tr a_{t} \,dt\big)^{1/2}\,dx<\infty\quad\text{(a.s.)},
\end{equation}
$$
  \int_{0}^{T}
\sup_{|x-x_{t}|\leq\gamma}\big(|G_{t}(x)|+|L_{t}F_{t}(x)|
+ 
|\Lambda_{t}F_{t}(x)|^{2}_{\ell_{2}}
$$
\begin{equation}
                                                           \label{4.6.4}
+|H_{t}(x)|^{2}_{\ell_{2}}+
\sum_{k=1}^{\infty}|
\Lambda^{k}_{t}H^{ k}_{t}(x)|
\big)\,dt<\infty\quad\text{(a.s.)}.
\end{equation}
 
Then for any $t\in[0,T]$ with probability one
$$
F _{t}(x_{t})=F _{0}(0)+\sum_{k=1}^{\infty}
\int_{0}^{t}(H^{ k}_{s}(x_{s})
+\Lambda^{k}_{s}F _{s}(x_{s}))\,dw^{k}_{s}
$$
\begin{equation}
                                          \label{3.31.1}
+\int_{0}^{t}\big(G _{s}(x_{s})
+L_{s}F _{s}(x_{s})+\sum_{k=1}^{\infty}
\Lambda^{k}_{s}H^{ k}_{s}(x_{s})\big)\,ds.
\end{equation}
Furthermore, if $F _{t}(x_{t})$ is continuous in $t$ 
on $[0,T]$ for almost
all $\omega$, then with probability one \eqref{3.31.1} holds
for all $t\in[0,T]$. 
\end{theorem}

Proof. First, observe that the series of stochastic integrals in \eqref{3.31.1}
converges uniformly on $[0,T]$
in probability and its limit is a local martingale
since
$$
\sum_{k=1}^{\infty}
\int_{0}^{T}|H^{ k}_{s}(x_{s})
+\Lambda^{k}_{s}F _{s}(x_{s})|^{2}\,ds
\leq2\sum_{k=1}^{\infty}
\int_{0}^{T}|H^{ k}_{s}(x_{s})|^{2}\,ds
$$
$$
+2\sum_{k=1}^{\infty}
\int_{0}^{T}|\Lambda^{k}_{s}F _{s}(x_{s})|^{2}\,ds
=2\int_{0}^{T}|H_{s}(x_{s})|^{2}_{\ell_{2}}\,ds
+2\int_{0}^{T} |\Lambda_{s}F_{s}(x_{s})|^{2}_{\ell_{2}}\,ds<\infty
$$
for almost all $\omega$ due to \eqref{4.6.4}. It is seen
 that, in light of \eqref{4.6.4}, the right-hand side
of \eqref{3.31.1} is continuous in $t$ (a.s.) and hence
the second assertion of the theorem follows from the first one.

To prove the first assertion, for $R\in(0,\infty)$ 
introduce  $\tau_{R}$ as the first exit time  of $x_{t}$
from $B_{R}$. Notice that if we take $x_{t\wedge\tau_{R}}$,
$\sigma _{t}I_{t<\tau_{R}}$, $b_{t}(x)I_{t<\tau_{R}}$,
$F_{t\wedge\tau_{R}}(x)$, $G_{t}(x)I_{t<\tau_{R}}$,
and $H_{t}(x)I_{t<\tau_{R}}$ instead of the original ones,
then the conditions (i)-(v) will be preserved. If we have 
\eqref{3.31.1} for the new objects, then we can send $R
\to\infty$ and easily obtain \eqref{3.31.1} as is 
 because $\tau_{R}\uparrow
\infty$, so that for  any $\omega$
there exists $R$ such that $t\wedge\tau_{R}=t$. 
We conclude that without loss of generality we may assume that,
for an  $R\in[0,\infty)$, we have
$|x_{t}|\leq R$ for all~$t$.

After that, by taking $\xi\in C^{\infty}_{0} $
such that it equals one on $B_{R+\gamma}$ and replacing
$F,G,H$ with $\xi F,\xi G,\xi H$, respectively, we see that
 the assumptions of the theorem will still be satisfied
and  assertion
\eqref{3.31.1} will be unaffected.
Therefore, without loss of generality we may assume that
there exists and $R<\infty$ such that
 $F,G,H$ vanish outside $B_{R}$.  
As in the proof of Lemma \ref{lemma 3.30.2} we may assume that
$\Omega'=\Omega$ in condition (iv). 

After these reductions we take a
 $\zeta\in C^{\infty}_{0} $, which is nonnegative,
radially symmetric,   with unit integral
and support in $B_{\gamma}$. 
Then for any $x\in\bR^{d}$   It\^o's formula yields that 
 with probability one
\begin{equation}
                                                 \label{3.20.1}
F_{t}(x)\zeta(x-x_{t})
-F_{0}(x)\zeta(x)=
\int_{0}^{t}\hat{G}_{s}(x)\,ds
+\sum_{k=1}^{\infty}\int_{0}^{t} \hat{H}^{k}_{s}(x)\,dw^{k}_{s}
\end{equation}
for all $t\in[0,T]$, where
$$
\hat{H}^{k}_{s}(x):=\zeta(x-x_{s})H^{k}_{s}(x)-
F_{s}(x)(\Lambda^{ k}_{s}
 \zeta)(x-x_{s}) ,\quad \hat{G}_{s}(x):=
\zeta(x-x_{s})G_{s}(x)
$$
$$
+F_{s}(x)(a^{ij}_{s}D_{ij}\zeta
-b^{i}_{s}D_{i}\zeta)(x-x_{s}) 
-\sum_{k=1}^{\infty}H^{k}_{s}(x)(\Lambda^{ k}_{s}
 \zeta) (x-x_{s}).
$$
We want to apply Lemma \ref{lemma 3.30.2} to
\eqref{3.20.1}.

First, observe that $\hat{G}$ and $\hat{H}$
satisfy assumption (i) of Lemma \ref{lemma 3.30.2} for any ball $\Gamma$
owing to the imposed measurability assumption 
on $G$ and $H$ and conditions (i) and (ii).
Then,  for each   $x\in\bR^{d}$ 
and $\omega\in\Omega$
$$
\int_{0}^{T}|\hat{G}_{t}(x)|\,dt
\leq N\int_{0}^{T}\eta_{t}(x)
\big( |G_{t}(x)|+|F_{t}(x)|(\tr a_{t}+|b_{t}|)\big)\,dt  
$$
$$
+N\int_{0}^{T}\eta_{t}(x) \sum_{k=1}^{\infty}|H^{k}_{t}(x)|\,
|\sigma^{k}_{t}|\,dt,
$$ 
where the constants $N$ are independent of $\omega, x$.
Regarding the last term notice that by  H\"older's inequality
$$
I(x):=\int_{0}^{T}\eta_{t}(x) \sum_{k=1}^{\infty}|H^{k}_{t}(x)|\,
|\sigma^{k}_{t}|\,dt
$$
$$
\leq
\big(\int_{0}^{T}\eta_{t}(x)|H _{t}(x)|^{2}_{\ell_{2}}
\,dt\big)^{1/2}\big(2\int_{0}^{T}\tr a_{t}\,dt\big)^{1/2},
$$
which is finite for all $\omega$ and $x$ owing to condition (iv)
and \eqref{3.30.7}. We see that by condition (iv)
$$
\int_{0}^{T}|\hat{G}_{t}(x)|\,dt<\infty
$$
for all $\omega,x$. Furthermore,
$$
I(x)\leq\big(\int_{0}^{T}
\sup_{y}(\eta_{t}(y)|H _{t}(y)|^{2}_{\ell_{2}})
\,dt\big)^{1/2}\big(2\int_{0}^{T}\tr a_{t}\,dt\big)^{1/2}
$$
which is finite (a.s.) due to assumption \eqref{4.6.4}.
The last expression here is independent of $x$
and the first one vanishes for $|x|\geq R$. Therefore, (a.s.)
$$
\int_{0}^{T}\int_{\bR^{d}}
\eta_{t}(x) \sum_{k=1}^{\infty}|H^{k}_{t}(x)|\,
|\sigma^{k}_{t}|\,dxdt<\infty.
$$
Similarly,  
$$
\int_{0}^{T}\int_{\bR^{d}}
\eta_{t}(x) |G_{t}(x)|\,dxdt<\infty
$$
(a.s.). By combining this with  
assumption \eqref{4.6.3}
  we see that  (a.s.)
$$
 \int_{0}^{T}\int_{\bR^{d}}|\hat{G}_{t}(x)|\,dxdt<\infty.
$$

Furthermore, condition (iv)
 implies that 
for each   $x\in\bR^{d}$ 
and $\omega\in\Omega$
$$
\int_{0}^{T}|\hat{H}_{t}(x)|_{\ell_{2}}^{2}\,dt<\infty
$$
and condition  (v) implies  that (a.s.)
$$
\int_{\bR^{d}}\big(\int_{0}^{T}|\hat{H}_{t}(x)|_{\ell_{2}}^{2}\,dt
\big)^{1/2}\,dx<\infty.
$$

We see that the assumptions of Lemma \ref{lemma 3.30.2}
are satisfied for any ball $\Gamma$ and recalling that 
 $F,G,H$ vanish for$|x|\geq R$ we conclude  
that for each $t\in[0,T]$ with probability one
$$
\int_{\bR^{d}} F_{t}(x)
\zeta(x-x_{t})\,dx=\int_{\bR^{d}} F_{0}(x)
\zeta(x)\,dx+\int_{0}^{t}\int_{\bR^{d}} 
\hat{G}_{s}(x)\,dxds
$$
$$
+\sum_{k=1}^{\infty}
\int_{0}^{t}\int_{\bR^{d}}\hat{H}^{k}_{s}(x)\,dxdw^{k}_{s}.
$$
 We fix $t\in[0,T]$ and 
use this formula with $\zeta_{\varepsilon}(x):
=\varepsilon^{-d}\zeta(x/\varepsilon)$, $\varepsilon>0$,
 in place of $\zeta$
and integrate by parts in the integrals  of $\hat{G}$ and
$\hat{H}$ with respect to $x$ (that is, use the definition
of generalized derivatives). Then by using the notation
$u^{(\varepsilon)}=u*\zeta_{\varepsilon}$ we find that
with probability one
$$
F^{(\varepsilon)}_{t}(x_{t})=F^{(\varepsilon)}_{0}(0)+
\sum_{k=1}^{\infty}
\int_{0}^{t}(H^{(\varepsilon)k}_{s}(x_{s})
+\Lambda^{k}_{s}F^{(\varepsilon)}_{s}(x_{s}))\,dw^{k}_{s}
$$
\begin{equation}
                                          \label{3.30.9}
+\int_{0}^{t}\big(G^{(\varepsilon)}_{s}(x_{s})
+L_{s}F^{(\varepsilon)}_{s}(x_{s})+\sum_{k=1}^{\infty}
\Lambda^{k}_{s}H^{(\varepsilon)k}_{s}(x_{s})\big)\,ds.
\end{equation} 

We now  let $\varepsilon\downarrow0$ in \eqref{3.30.9}.
Since $F_{t}(x)$ is continuous in $x$ we have
$F^{(\varepsilon)}_{t}(x_{t})\to F _{t}(x_{t})$
(for all $\omega$).
Furthermore, by assumption for almost any $\omega$, for almost all
$s\in[0,t]$ the function $|\Lambda_{s}F_{s}(x)|_{\ell_{2}}$
is continuous and also $\Lambda^{k}_{s}F_{s}(x)$ are continuous.
It follows by Dini's theorem that for the above $\omega$ and $s$
\begin{equation}
                                                           \label{4.7.1}
 \sum_{k=n}^{\infty}|\Lambda^{k}_{s}F_{s}(x)|^{2}
\downarrow0
\end{equation}
as $n\to\infty$ uniformly on compact sets in $\bR^{d}$.
 Similar argument shows that
(a.s.) for almost any $s\in[0,t]$
$$
\sum_{k=n}^{\infty}(| 
 H^{k}_{s}| ^{2}+|
\Lambda^{k}_{s} H^{ k}_{s}|)(x)\to0,
$$
as $n\to\infty$ uniformly on compact sets in $\bR^{d}$.
This implies that (a.s.) for almost any $s\in[0,t]$ as 
$\varepsilon\downarrow0$,
$$
\big(|(\Lambda_{s}F_{s} )^{(\varepsilon)}
- \Lambda_{s}F_{s} |_{\ell_{2}} ^{2} 
+|H^{(\varepsilon)}_{t}-H _{t}|^{2}_{\ell_{2}}
+\sum_{k=1}^{\infty}|(\Lambda^{k}_{s}  H^{(\varepsilon) k}_{s}
-\Lambda^{k}_{s} H^{ k}_{s}|\big)(x_{s})
\to0.
$$

Hence,
in light of \eqref{4.6.4},
by the dominated convergence theorem (a.s.)
  we have as $\varepsilon\downarrow0$ that
$$
\int_{0}^{t}\big(|(\Lambda_{s}F_{s} )^{(\varepsilon)}
- \Lambda_{s}F_{s} |_{\ell_{2}} ^{2} 
+|H^{(\varepsilon)}_{t}-H _{t}|^{2}_{\ell_{2}}
+\sum_{k=1}^{\infty}| \Lambda^{k}_{s} H^{(\varepsilon) k}_{s} 
-\Lambda^{k}_{s} H^{ k}_{s}|\big)(x_{s})\,ds 
\to0.
$$
  This allows us to assert
that part of the terms in \eqref{3.30.9}
converges in probability to what we need. 

Convergence of the remaining terms in \eqref{3.30.9}
is proved in like manner. Thus, passing to
the limit in \eqref{3.30.9} yields \eqref{3.31.1}
and this brings the proof of the theorem to an end.

\begin{remark}
The assumptions of this theorem are substantially weaker
than the ones usually imposed (see, for instance, Theorem 3.3.1
of \cite{Ku} or Theorem 1.4.9 of \cite{Ro}). In particular,
we are dealing with the generalized functions
$L_{t}F_{t}(x)$, $\Lambda^{k}_{t}F_{t}(x)$, and
$\Lambda^{k}_{t}H^{k}_{t}(x)$, which always exist
and are continuous in $x$ (just equal zero) at those
$(t,\omega)$ at which   $a_{t}=b_{t}=0$.
Furthermore, at those
$(t,\omega)$, at which  $a_{t}=0$,  no differentiability assumption
is imposed on $H$.

Also, observe that if 
  $G\equiv0, H\equiv0$, then the It\^o-Wentzell
formula becomes just It\^o's formula and our theorem gives the proof
of it under substantially weaker assumptions
than the ones usually imposed.
\end{remark}

\mysection{Proof of Theorem \protect\ref{theorem 11.16.5}}
                                           \label{section 4.6.1}

 Here we suppose that the assumptions of Theorem 
\ref{theorem 11.16.5} are satisfied
  and will base our proof on 
Theorem \ref{theorem 3.30.1}.

Note that
\eqref{4.4.5} involves the values of $x_{t}$ only
for $t<\tau$. Therefore, without losing generality
we may assume that $b^{i}_{t}=\sigma^{ik}_{t}=0$
for $t\geq\tau$ for all $i,k$. Also obviously we may assume 
that $\tau$ is bounded. 
These additional assumptions
are supposed to hold throughout the section.
 
\begin{lemma}
                                          \label{lemma 4.10.1}
 Take a $\phi\in C^{\infty}_{0} $ and set
$$
F_{t}(x) =(u_{t\wedge\tau}( \cdot+x),\phi).
$$
Then with probability one $F_{t}(x)$  
and, for any multi-index $\alpha$, $D^{\alpha}F_{t}(x)$
are  a continuous functions
of $(t,x)\in\bR_{+}\times\bR^{d}$.
\end{lemma}

Proof. 
By multiplying, if necessary, $u,f,g$ by the indicator
of an event of full probability (perhaps depending on $\phi$)
 we may assume that
for any  $\omega\in\Omega$, multi-index $\alpha$, and 
 $T, R\in\bR_{+}$, we have
$$
\int_{0}^{T}\sup_{ |x|\leq R} |(f_{t} ,
D^{\alpha}\phi(\cdot-x))|I_{t\leq\tau} \,dt<\infty,
$$
$$
 \int_{0}^{T}\sup_{ |x|\leq R} |(g_{t} ,
D^{\alpha}\phi(\cdot-x))|_{\ell_{2}}^{2}I_{t\leq\tau}\,dt<\infty.
$$ 
Set 
$$
 G_{t}(x)
=(f_{t}( \cdot+x),\phi)I_{t\leq\tau},
\quad
 H _{t}(x)=(g _{t}( \cdot+x),\phi)I_{t\leq\tau}
$$
and observe that   for each $x\in\bR^{d}$
equation \eqref{4.9.1} holds   with 
probability one for all $t$ due to \eqref{12.23.40}. 
Therefore,
\begin{equation}
                                                        \label{4.12.4}
m_{t}(x):=F_{t}(x)-\int_{0}^{t}G_{s}(x)\,ds
\end{equation}
for every $x$ with probability one satisfies
$$
m_{t}(x)=\sum_{k=1}^{\infty}\int_{0}^{t}H^{k}_{s}(x)\,dw^{k}_{s}
$$
for all $t$. Furthermore, for each $\omega$ and $t$ the functions
$F_{t}(x)$, $G_{t}(x)$, and $H^{k}_{t}(x)$ are infinitely differentiable
in $x$.   If $|x|\leq R$ and $\alpha$ is a multi-index,
 then
$$
 |D^{\alpha}G_{s}(x)|\leq\sup_{|y|\leq R}|(f_{s}(\cdot+y),D^{\alpha}\phi)|
I_{s\leq\tau}
$$
and by assumption the latter is locally integrable on $\bR_{+}$
(for any $\omega$).
This   shows that the integral
in \eqref{4.12.4} is also infinitely differentiable
in $x$ and one can perform  differentiating the integral
by differentiating  the integrand. Hence $m_{t}(x)$ is infinitely 
differentiable in $x$. By replacing $\phi$ with $D^{\alpha}\phi$
in the above argument  we now see that for every $x$ 
 and any multi-index $\alpha$
with probability one  
$$
D^{\alpha}m_{t}(x)=
\sum_{k=1}^{\infty}\int_{0}^{t}D^{\alpha}H^{k}_{s}(x)\,dw^{k}_{s}
$$
for all $t$. 
The quadratic variation of the sum on the right can be
taken to be
$$
A^{\alpha}_{t}(x)=\sum_{k=1}^{\infty}\int_{0}^{t}I_{s\leq\tau}
( g ^{k}_{s}, D^{\alpha}\phi(\cdot-x))^{2}\,ds.
$$
Moreover,
 for each $x\in\bR^{d}$ the function $D^{\alpha}m_{t}(x)$ is
$\cF\otimes\cB(\bR_{+})$-measurable (see Definition \ref{def 10.25.1})
and for each $(\omega,t)$ continuous in $x$. Hence it is
a regular martingale field on $\bR^{d}$. Finally,
by definition for any ball $\Gamma$ (recall that $\tau$ is bounded)
$$
\big(\int_{\Gamma}|A^{\alpha}_{\infty}(x)|^{1/2}\,dx\big)^{2}
\leq|\Gamma|^{2}\sup_{x\in \Gamma }\int_{0}^{\tau}
|( g _{s},D^{\alpha}\phi(\cdot-x))|_{\ell_{2}}^{2}\,ds
$$
$$
\leq|\Gamma|^{2}\int_{0}^{\tau}\sup_{x\in \Gamma }
|(g_{s},D^{\alpha}\phi(\cdot-x))|_{\ell_{2}}^{2}\,ds<\infty,
$$
where $|\Gamma|$ is the volume of $\Gamma$.
Now our assertion  about $F$  follows immediately
from Corollary \ref{corollary 4.3.1}. 
 The functions  $D^{\alpha}F$ are taken care of
by replacing $\phi$ with $D^{\alpha}\phi$.  
The lemma is proved.

\begin{remark}
                                           \label{remark 4.3.2}
Naturally, equation \eqref{4.4.5} is understood in the sense
of Definition \ref{def 10.25.3}. Therefore, it is important to
explain that the terms on the right
 in \eqref{4.4.5} belong to the right class
of functions.
Notice that
 for any $\phi\in C^{\infty}_{0}$    
  and $T,R\in\bR_{+}$ we have
$$
\int_{0}^{T}\sup_{|x|\leq R}
|(D_{1}v_{t}(\cdot+x),\phi)\sigma^{1\cdot}_{t}|_{\ell_{2}}^{2}\,dt
=
\int_{0}^{T}\sup_{|x|\leq R}|(v_{t}(\cdot+x), D_{1} \phi)|^{2} a^{11}_{t}\,dt
$$
$$
\leq  \sup_{|x|\leq R,t\leq T}|
(u_{t\wedge\tau}(\cdot+x_{t\wedge\tau}+x, D_{1} \phi)|^{2} 
\int_{0}^{\tau}a^{11}_{t}\,dt
$$
$$
\leq \sup_{|x|\leq R+N,t\leq T}|
(u_{t\wedge\tau}(\cdot +x, D_{1} \phi)|^{2} 
\int_{0}^{\tau}a^{11}_{t}\,dt<\infty\quad\hbox{(a.s.)},
$$
where $N=\sup_{t} |x_{t\wedge\tau}|$ is a finite random
variable and the last inequality follows from Lemma \ref{lemma 4.10.1}
and assumption \eqref{3.30.7}.
Similarly one treats the remaining  terms in \eqref{4.4.5}.

\end{remark}

{\bf Proof of Theorem \ref{theorem 11.16.5}}.  We take $\phi,F,G$, and $H$
from Lemma \ref{lemma 4.10.1}.
By definition, 
$G_{t}(x)$  and $H^{k}_{t}(x)$
   are predictable
for each $x$. In addition, for each   $\omega$  and $t$ these functions 
are infinitely 
differentiable with respect to $x$ (in the strong sense in $\ell_{2}$
in the case of $H_{t}(x)$). Therefore,
these functions satisfy
the measurability condition (i) of Lemma \ref{lemma 3.30.2} for any $T$ and
ball $\Gamma$  and satisfy condition (iii) (a)
of Theorem \ref{theorem 3.30.1}. Similarly   $F_{t}(x)$ satisfies conditions
  (i) and (ii) of Theorem \ref{theorem 3.30.1} for any $T\in\bR_{+}$.

Furthermore,   not only  $H^{k}_{t}(x)$ and
 $|H_{t}(x)|_{\ell_{2}}$ are continuous, the same is true for
the derivatives of $H_{t}(x)$ with respect to $x$. In particular,
as in the case of \eqref{4.7.1}, for $i=1,2,...,d$
$$
\sum_{k=n}^{\infty}|D_{i}H^{k}_{t}(x)|^{2}\downarrow0
$$
as $n\to\infty$ uniformly on compact subsets of $\bR^{d}$
for any $\omega$ and $t$. The estimate
\begin{equation}
                                                    \label{4.7.3}
 \sum_{k=n}^{\infty}|
\Lambda^{k}_{t}H^{k}_{t}(x)|\leq(2\tr a_{t})^{1/2}
 \big(\sum_{i=1}^{d}\sum_{k=n}^{\infty}
|D_{i}H^{k}_{t}(x)|^{2}\big)^{1/2}
\end{equation}
obtained by   H\"older's inequality, shows that the left-hand side
goes to zero as $n\to\infty$ uniformly on compact subsets of $\bR^{d}$
 and hence its value for $n=1$
is a continuous function of $x$  
for   any $\omega$ and $t$. Also, the mean value theorem
and the continuity of the second order derivatives of $F_{t}(x)$
in $x$ easily yield the continuity of 
$|\Lambda_{t}F_{t}(x)|_{\ell_{2}}$.

Thus, assumptions (ii) and (iii) of Theorem \ref{theorem 3.30.1} 
are satisfied for any $T\in\bR_{+}$.
In assumption (iv) equation \eqref{4.9.1} holds for each $x$ with 
probability one for all $t$ due to \eqref{12.23.40}.

Furthermore, by assumption for each $T, R \in\bR_{+}$, $\omega$,
 and $|x|\leq R$  we have
\begin{equation}
                                                    \label{4.4.3} 
 \int_{0}^{T}|G_{t}(x)|\,dt 
\leq\int_{0}^{T}\sup_{ |y|\leq R }|G_{t}(y)|\,dt<\infty.
\end{equation}
Inequalities like \eqref{4.7.3} and \eqref{4.4.3},
 Lemma \ref{lemma 4.10.1},  and the fact that
the support of $\sup_{t\leq T}\eta_{t}(x)$
is bounded for  
 each $\omega$ imply that that assumptions
(iv) and (v) of Theorem \ref{theorem 3.30.1} are satisfied
and we can apply it.

 Now by Lemma \ref{lemma 4.10.1} and Theorem \ref{theorem 3.30.1}  
we have that with
probability one  \eqref{3.31.1} holds for all $t\in[0,T]$
and, actually, by the arbitrariness of $T$, for all $t$.
We rewrite \eqref{3.31.1} in terms of $u_{t}$ and $v_{t}$,
use that $b^{i}_{t}=\sigma^{ik}_{t}=0$
for $t\geq\tau$ for all $i,k$, 
and see that with probability one
$$
(v_{t\wedge\tau},\phi)=(v_{0},\phi)
+\sum_{k=1}^{\infty}
\int_{0}^{t }I_{s\leq\tau}( (g^{k}_{s}
( \cdot + x_{s}),\phi)
 + (\sigma^{ik}_{s}D_{i}v_{s},\phi))\,dw^{k}_{s}
$$
\begin{equation}
                                                \label{4.4.6}
+\int_{0}^{t }I_{s\leq\tau}
\big[(f_{s}( \cdot + x_{s})+ L_{s} v_{s},\phi)
 + \sum_{k=1}^{\infty}\sigma^{ik}_{s}(D_{i}g^{k}_{s}(\cdot
 + x_{s}),\phi)\big]\,ds
\end{equation}
for all $t$. Here for each $\omega$ and $s$     
(recall the definition of the limit of distributions) 
$$
\sum_{k=1}^{\infty}\sigma^{ik}_{s}(D_{i}g^{k}_{s}(\cdot
 + x_{s}),\phi)=
\big(\sum_{k=1}^{\infty}\sigma^{ik}_{s} D_{i}g^{k}_{s}(\cdot
 + x_{s}),\phi\big)
$$
since
$$
\big(\sum_{k=1}^{\infty}|\sigma^{ik}_{s}|\,|(D_{i}g^{k}_{s}(\cdot
 + x_{s}),\phi)|\big)^{2}\leq2\tr a_{s}\sum_{i=1}^{d}
|(D_{i}g_{s}(\cdot
 + x_{s}),\phi)|^{2}_{\ell_{2}}<\infty.
$$
This shows that \eqref{4.4.6}  implies \eqref{4.4.5}.
The theorem is proved.

\end{document}